\documentclass{amsart2000}
\input diagrams
\diagramstyle[PostScript=dvips]

\newcommand{\nc}{{\mathbb{C}}}        
\newcommand{\nz}{{\mathbb{Z}}}        

\newcommand{\tensor}{\otimes}
\newcommand{\cross}{\times}

\newcommand{\M}{\overline{\MM}}       
\newcommand{\MM}{\mathcal{M}}          

\newcommand{\mgnbar}{\M_{g,n}}
\newcommand{\bg}{\cb G}               

\newcommand{\cb}{{\mathcal B}}

\newcommand{\ch}{{\mathcal H}}

\newcommand{\cx}{{\mathcal X}}
\newcommand{\cy}{{\mathcal Y}}

\newcommand{\bgamma}{\boldsymbol{\gamma}}
\newcommand{\bmu}{\boldsymbol{\mu}}

\newcommand{\sigmas}{{\sigma_1, \sigma_2, \sigma_3}}
\newcommand{\bga}{\bgamma}
\newcommand{\bsig}{\boldsymbol{\sigma}}

\newcommand{\ad}{\operatorname{ad}}
\newcommand{\class}{\operatorname{Class}}
\newcommand{\cft}{CohFT}              
\newcommand{\cfts}{CohFTs}            

\newcommand{\cor}[1]{\langle {#1} \rangle}                 
\newcommand{\ccor}[1]{\langle\langle {#1} \rangle\rangle}  

\newcommand{\bgt}{\widetilde{\bg}}

\newcommand{\bt}{\mathbf{t}}
\newcommand{\btau}{\boldsymbol{\tau}}
\newcommand{\bu}{\mathbf{u}}
\newcommand{\corG}[1]{\left< {#1} \right>^G} 
\newcommand{\PhiG}{\Phi^G}
\newcommand{\tbu}{\widetilde{\bu}}
\newcommand{\tu}{\widetilde{u}}
\newcommand{\ZG}{Z^G}

\newtheorem{thm}{Theorem}[section]
\newtheorem{lm}[thm]{Lemma}
\newtheorem{prop}[thm]{Proposition}
\newtheorem{crl}[thm]{Corollary}

\theoremstyle{definition}
\newtheorem{rem}[thm]{Remark}

\theoremstyle{remark}
 
\newtheorem{ack}{Acknowledgments\hspace{-.5em}}  

\numberwithin{equation}{section}

\begin{document}

\copyrightinfo{2002}{Tyler J. Jarvis and Takashi Kimura}

\title[Orbifold quantum cohomology of $\bg$]
{Orbifold quantum cohomology of the classifying space of a finite
group}

\subjclass[2000] {Primary: 14N35, 53D45. Secondary: 20C05, 14H10}

\author
[T. J. Jarvis]{Tyler J. Jarvis}
\address
{Department of Mathematics, Brigham Young University, Provo, UT
84602, USA}
\email{jarvis@math.byu.edu}
\thanks{Research of the first author was partially supported by NSF
grant DMS-0105788}

\author
[T. Kimura]{Takashi Kimura}
\address
{Department of Mathematics, 111 Cummington Street, Boston
University, Boston, MA 02215, USA} \email{kimura@math.bu.edu}

\date{}

\begin{abstract}
We work through, in detail, the quantum cohomology, with gravitational descendants,
of the orbifold $\bg$, the point with action of a finite group $G$. We provide a
simple description of algebraic structures on the state space of
this theory.  As a consequence, we find that multiple copies of
commuting Virasoro algebras appear which completely determine the correlators
of the theory.

\end{abstract}

\maketitle

\section{Introduction}
\label{intro}

W. Chen and Y. Ruan \cite{CR1, CR2} introduced the notion of the
Gromov-Witten invariants of an orbifold $V$. Their construction
reduces to the usual  Gromov-Witten invariants when $V$ is a
smooth, projective variety. When $V = [Y/G]$, where $Y$ is a smooth,
projective variety, and $G$ is a finite group, the state space of
this theory is generally larger than the $G$-invariant part of the
cohomology of $Y$; indeed, it has additional direct summands
associated to loci in $Y$ with nontrivial isotropy. These loci are
called the twisted sectors of the theory, and their presence
should be part of the proper notion of the cohomology of an
orbifold.

The purpose of this paper is to provide a detailed treatment of
this theory for the simplest case, namely, when $V$ is the
classifying stack $\bg$ of a finite group $G$. Many of the features and
subtleties of the Gromov-Witten invariants of orbifolds are
present even here. The state space of this theory contains twisted
sectors and the correlators in this theory are intersection
numbers on $\M_{g,n}(\bg)$,  the moduli space of genus-$g$,
$n$-pointed orbifold stable maps into $\bg$.

The correlators in this theory can be described in purely group
theoretic terms, and we recover the result \cite{Ru} that the
algebraic structure on the state space $\ch$ is isomorphic to the
center of the group algebra,  $Z\nc[G]$, together with an
invariant metric.

Furthermore, we show that on the large phase space of this theory,
there are $r$ commuting copies of ``half'' the Virasoro algebra in
this theory, where $r$ is the dimension of $\ch$, all of which
annihilate, and completely determine, the exponential of the large
phase space potential function. We obtain a proof of the usual
Virasoro conjecture \cite{EHX} as the special case where a
diagonal action (after a variable rescaling) of these Virasoro
algebras is considered.  Similarly, the relevant integrable
hierarchy consists of $r$ commuting copies of the KdV hierarchy.

Finally, it is worth pointing out that the moduli spaces
$\M_{g,n}(\bg)$ have many features similar to $\M_{g,n}^{1/r}$,
the moduli space of $r$-spin curves \cite{JKV}. The moduli spaces
$\M_{g,n}(\bg)$ have boundary strata indexed by stable graphs
whose tails and half-edges are decorated by elements of $G$ (up to
conjugation) while $\M_{g,n}^{1/r}$ have boundary strata indexed
by stable graphs whose tails and half-edges are decorated by
elements of $\nz_{r}$. The construction of the correlators in both
theories is also analogous. However, the ring structures on the
state spaces are distinctly different.  It would be interesting to
find the analogs of the above results for the spaces
$\M_{g,n}^{1/r}$.

\subsection{Conventions and notation}

Throughout this paper, except where otherwise specified, we will
work only over the complex numbers $\nc$, and all groups will be
finite.

For a given group $G$ we denote the classifying stack of $G$  by
$\bg$; namely, $\bg$ is the stack quotient $$\bg := [pt/G]$$ of a
point modulo a trivial $G$ action.

Group elements will always be denoted by lower-case Greek letters,
and the conjugacy class of an element $\gamma \in G$ is denoted
$[\![\gamma]\!]$.  An $n$-tuple of elements
$(\gamma_1,\dots,\gamma_n)$ will generally be denoted by a
boldfaced $\bga$.  The centralizer in $G$ of an element $\gamma$
is denoted $C(\gamma)$, and the intersection $\bigcap_{i =1}^n
C(\gamma_i)$ of several centralizers is denoted
$C(\gamma_1,\dots,\gamma_n) = C(\bga)$. The commutator
$\alpha\beta\alpha^{-1}\beta^{-1}$ of two group elements $\alpha$
and $\beta$ is denoted $[\alpha,\beta]$. Finally, the center of an
algebra $A$ will be denoted $ZA$.

\begin{ack}
We would like to thank Dan Abramovich and Yongbin Ruan for many
helpful discussions, and Heidi Jarvis for help with typesetting.
\end{ack}

\section{Orbifold stable maps into $\bg$}

\subsection{The stack}

Our chief objects of study in this paper are Gromov-Witten
invariants arising from the stack $\mgnbar(\bg)$ of $n$-pointed
orbifold stable maps into $\bg$, in the sense of W. Chen and Y.
Ruan \cite{CR1, CR2} (called balanced twisted stable maps by
Abramovich and Vistoli \cite{AV}).

These are maps $f:\Sigma \rTo \bg$ from an orbifold Riemann
surface (orbicurve) $\Sigma$ into the classifying stack $\bg$ of a
finite group $G$.  Here $\Sigma$ has non-trivial orbifold
structure only at marked points $p_1, \dots, p_n$ and at nodes,
and the orbifold structure at the nodes is \emph{balanced},
meaning that the  action of the stabilizer $G_q \simeq \bmu_l$ at
a nodal geometric point $q$ of $\Sigma$ has complementary
eigenvalues on the tangent spaces of the two branches of $\Sigma$
at $q$.

The stack $\mgnbar(\bg)$ is a smooth, proper Deligne-Mumford stack
with projective coarse moduli space \cite[Thm 3.0.2]{ACV}.  This
stack has a number of important connections to other moduli
problems.  For example, in the case that $G$ is the symmetric
group $S_d$ on $d$ letters, the stack $\mgnbar(\cb S_d)$ is the
normalization of the stack of admissible covers \cite[\S4]{ACV}.

Recall that an orbifold stable map from a smooth, $n$-pointed
orbicurve $(\Sigma, p_1,\dots, p_n)$ into $\bg$ determines a
principal $G$ bundle on the complement $\Sigma-\{p_1, \dots,
p_n\}$.  This, in turn, is determined by its holonomy, that is by
a homomorphism $\pi_1(\Sigma-\{p_1, \dots, p_n\}) \rTo G$.
Moreover, the stabilizer $G_{p_i}$ of a marked point $p_i$ of
$\Sigma$ is always cyclic, and the order of $G_{p_i}$ is equal to
the order of the holonomy around that marked point.  Conversely,
$G$ acts by conjugation on the homomorphisms $\pi_1(\Sigma-\{p_1,
\dots, p_n\}) \rTo G$, and two such homomorphisms determine the
same $G$-bundle on $\Sigma-\{p_1, \dots, p_n\}$ precisely when
they differ by this adjoint action  of $G$.

Finally, since $\bg$ is a proper, separated stack, any principal
$G$ bundle on $\Sigma-\{p_i,\dots,p_n\}$  extends uniquely, after
suitable base extension, to a principal bundle on some proper
curve $\tilde{\Sigma}$, with an isomorphism
$\tilde{\Sigma}-\{p_1,\dots,p_n\} \rTo \Sigma-\{p_1,\dots,p_n\}$.
The data of this cover (up to the obvious notion of equivalence of
such data) is exactly equivalent to the data of a principal bundle
on an orbicurve.   Moreover, the adjoint action of $G$ on these
homomorphisms exactly corresponds to the natural action of $G$ on
the orbifold stable maps.   Thus we have the following:

\begin{prop}\label{point}
For a given curve $[C, p_1,\dots,p_n]$ corresponding to a point of
the smooth locus $\MM_{g,n}$, the fiber of $\mgnbar(\bg)$ over the
point $[C]$ corresponds to the quotient
$\operatorname{Hom}(\pi_1(C-\{p_1,\dots,p_n\}, G)/ \ad G$.
\end{prop}

\subsection{Morphisms}

There are several natural morphisms from $\mgnbar(\bg)$.

First, recall that for $X=\bg$ the twisted sectors
$X_{[\![\gamma]\!]} = \{ (x, [\![\gamma]\!]) | x \in X, \gamma\in
G_x\}$ of \cite{CR1} are simply  $\bg_{[\![\gamma]\!]} \cong
[pt/C(\gamma)]$.

There are evaluation morphisms $$ev_i: \mgnbar(\bg) \rTo \bgt,$$
where ${\bgt} = \coprod_{[\![\gamma]\!]} \bg_{[\![\gamma]\!]}$ is
the disjoint union of twisted sectors.  We can describe the
evaluation morphism $ev_i$ as follows:  Any stable map $f:\Sigma
\rTo \bg$ must be representable, and hence must induce an
injective homomorphism $f_*:G_{p_i} \rTo G$ from the local group
$G_{p_i}=\nz/m_i$ of the $i$th marked point $p_i$ of $\Sigma$ into
$G$.  The group $G$ acts by conjugation on this homomorphism
$f_*$, and so the image of $1\in \nz/m_i$ in $G$ is defined only
up to conjugacy.  The evaluation morphism $ev_i$ is the morphism
taking $[\Sigma,p_1,\dots,p_n, f]$ to the point $(f(p_i),
[\![f_*(1)]\!]) \in \bg_{[\![f_*(1)]\!]} \subseteq {\bgt}$.
Alternatively, the image $f_*(1)$ is simply the holonomy of the
induced $G$-bundle on $\Sigma-\{p_1, \dots,p_n\}$ around the
marked point $p_i$.

We can use the evaluation morphism to see that the stack
$\M_{g,n}(\bg)$ breaks up as the disjoint union of open and closed
substacks $$\M_{g,n}(\bg) =
\coprod_{([\![\gamma_1]\!],\dots,[\![\gamma_n]\!])}
\M_{g,n}(\bg,[\![\gamma_1]\!],\dots,[\![\gamma_n]\!]),$$ where
$\M_{g,n}(\bg,[\![\gamma_1]\!],\dots,[\![\gamma_n]\!])=
ev_1^{-1}(\bg_{[\![\gamma_1 ]\!]}) \cap \dots \cap
ev_n^{-1}(\bg_{[\![\gamma_n ]\!]})$ is the substack of points of
$\M_{g,n}(\bg)$ mapped by $ev_i$ to $\bg_{\gamma_i}$ for every
$i\in\{1,\dots,n\}$. Of course,
$\M_{g,n}(\bg,[\![\gamma_1]\!],\dots,[\![\gamma_n]\!])$ may be empty for some
choices of conjugacy classes $([\![\gamma_1]\!],\dots,[\![\gamma_n]\!])$.

In the special case of $\M_{0,3}(\bg)$, since there is only one
$3$-pointed, genus-zero stable curve (call it $\Sigma$), we may
fix, once and for all a base point $q$ (distinct from the three
marked points $p_1$, $p_2$, and $p_3$) and a basis $\{s_1, s_2,
s_3\}$ of $\pi_1(\Sigma-\{p_1,p_2,p_3\}, q)$ such that
$\prod_{i=1}^3 s_i = 1$.   In this case, Proposition~\ref{point}
shows that each component of $\M_{0,3}(\bg)$ is uniquely
determined up to simultaneous (diagonal) adjoint action of $G$ by
a choice of three elements $\gamma_1, \gamma_2, \gamma_3 \in G$
such that $\prod \gamma_i = 1$.  That is, if we let
$[\![\bgamma]\!]$ denote the diagonal conjugacy class of the
triple $\bgamma = (\gamma_1, \gamma_2,\gamma_3)$, and we let $T^3$
denotes the set of all such triple conjugacy classes, whose
product is trivial (i.e., $\prod \gamma_i = 1$), then
$$\M_{0,3}(\bg) = \coprod_{[\![\bgamma]\!]\in T^3}
\M_{0,3}(\bg,[\![\bgamma]\!]),$$ and for a given
$[\![\bgamma]\!]\in T^3$, we have $$\M_{0,3}(\bg,[\![\bgamma]\!])
\cong \mathcal{B}C(\bgamma).$$

For $g>0$ or $n>3$ one cannot index the moduli space so easily,
since deformations in the moduli space may act on the holonomies
in non-diagonal ways.  Nevertheless, for any fixed, smooth,
$n$-pointed curve, Proposition~\ref{point} gives a complete
description of the points of $\M_{g,n}(\bg)$ that lie over it.
This shows that the forgetful morphism   $\pi:\mgnbar(\bg) \rTo
\mgnbar$, is quasi-finite.  In fact $\pi$ is also proper, but it
is not generally representable \cite[Cor 3.0.5]{ACV}.

Finally we have the  forgetting-tails morphism, defined as
follows.  When the holonomy $\gamma_i$  around a marked point
$p_i$ is trivial, then we may \emph{forget} the data of that
marked point.  This gives a morphism $$\M_{g,n+1}(\bg,
[\![\gamma_1]\!], \dots, [\![1]\!], \dots, [\![\gamma_n]\!]) \rTo
\mgnbar(\bg, [\![\gamma_1]\!],\dots, \widehat{[\![1]\!]}, \dots,
[\![\gamma_n]\!]).$$
Note that the forgetting-tails morphism is
not defined for all of $\M_{g,n+1}(\bg)$, but rather only for the
components corresponding to marked points with trivial holonomy.

\section{Gromov-Witten invariants, cohomological field theory and K-theory}

We define the classes $\psi_i$ on $\mgnbar(\bg)$ to be the
pullbacks $\psi_i = \pi^*(\psi_i)$ of the $\psi_i$ classes on
$\mgnbar$.

The tangent bundle of $\bg$ is trivial, and thus the virtual
fundamental class of $\mgnbar(\bg)$ is just the usual fundamental
class.

The orbifold cohomology of $\bg$ is, as a vector space, $$\ch :=
H^*_{orb}(\bg,\nc) := H^*({\bgt,\nc}) =
\bigoplus_{[\![\gamma]\!]} \nc.$$

For each conjugacy class $[\![\gamma]\!]$ in $G$, let
$e_{[\![\gamma]\!]}$ denote the class $1 \in
H^0(\bg_{[\![\gamma]\!]}, \nc) \subseteq \ch$. The
$e_{[\![\gamma]\!]}$'s form a basis of $\ch$ and we may form
\emph{$n$-point correlators}
$$\cor{\tau_{a_1}(e_{[\![\gamma_1]\!]})\dots
\tau_{a_n}(e_{[\![\gamma_n]\!]})}^G_g := \int_{\mgnbar(\bg)}
\prod_i \psi_i^{a_i} ev_i^*(e_{[\![\gamma_i]\!]}).$$

\subsection{Cohomological field theory}

The three-point, genus-zero correlators play a special role, since
they define the metric on $\ch$ and the quantum (orbifold)
product.  They vanish for dimensional reasons unless $\sum_{i=1}^n
a_i =0$.

\begin{prop}
We have \begin{align*} \langle\tau_0(e_{[\![\gamma_1]\!]})
\tau_0(e_{[\![\gamma_2]\!]})
\tau_0(e_{[\![\gamma_3]\!]})\rangle^G_0 & =
\sum_{\substack{\sigmas \\ \prod \sigma_i=1 \\ \sigma_i \in
[\![\gamma_i]\!]}} \frac{1}{|G|}\\ & =
\sum_{\substack{[\![\sigma_1, \sigma_2, \sigma_3]\!] \\ \prod
\sigma_i =1 \\ \sigma_i \in [\![\gamma_i]\!]}}
\frac{1}{|C(\sigma_1, \sigma_2)|}. \end{align*}
\end{prop}

\begin{proof}
For a given component $\M_{0,3} (\bg, [\![\sigma_1, \sigma_2,
\sigma_3]\!])$ of the moduli space $\M_{0,3}(\bg)$, the evaluation
map $ev_i$ maps this component to $\bg_{[\![\sigma_i]\!]}$, so
$$ev_i^*(e_{[\![\gamma_i]\!]})= \left\{ \begin{array}{ccc} 0 &
\text{if} & \sigma_i \notin [\![\gamma_i]\!]
\\ 1 & \text{ if } & \sigma_i \in [\![\gamma_i]\!]\end{array} \right. .$$
Moreover, $$\int_{\M_{0,3}(\bg,[\![\sigma_1, \sigma_2,
\sigma_3]\!])} 1 = \deg(\M_{0,3}(\bg,[\![\sigma_1, \sigma_2,
\sigma_3]\!]))=\frac{1}{|C(\sigmas)|} = \frac{1}{|C(\sigma_1,
\sigma_2)|}.$$  So
\begin{equation}\label{tau} \langle\tau_0(e_{[\![\sigma_1]\!]})
\tau_0(e_{[\![\sigma_2]\!]})
\tau_0(e_{[\![\sigma_3]\!]})\rangle^G_0 =
\sum_{\substack{[\![\sigmas]\!]\\ \prod \sigma_i=1 \\ \sigma_i \in
[\![\gamma_i]\!]}} \frac{1}{|C(\sigma_1, \sigma_2)|},
\end{equation} and it is easy to see that the number of elements
in a given non-empty conjugacy class $[\![\sigmas]\!]$ of triples
is exactly $\frac{|G|}{|C(\sigma_1, \sigma_2)|}$, so the
expression (\ref{tau}) gives the rest of the proposition.
\end{proof}

\begin{crl}
The metric on $\ch$ induced by the 3-point correlators is
$$\eta_{[\![\gamma_1]\!][\![\gamma_2]\!]} :=
\eta(e_{[\![\gamma_1]\!]},e_{[\![\gamma_2]\!]}) =
\frac{1}{|C(\gamma_1)|}
\delta_{[\![\gamma_1]\!][\![\gamma^{-1}_2]\!]},$$ which is
non-degenerate on $\ch$.

The quantum product is given by $$e_{[\![\gamma_1]\!]} \ast
e_{[\![\gamma_2]\!]}= \sum_{\substack{\sigma_1, \sigma_2 \\
\sigma_i \in [\![\gamma_i]\!]}} \frac{|C(\sigma_1 \sigma_2)|}{|G|}
e_{[\![\sigma_1 \sigma_2]\!]}.$$
\end{crl}

It is clear that $\ch$ is additively isomorphic to the ring
$\class_{\nc}(G)$ of class functions of $G$, where
$e_{[\![\gamma]\!]}$ is the class function
$e_{[\![\gamma]\!]}(\sigma) =
\delta_{[\![\gamma]\!],[\![\sigma]\!]}$.  Moreover,
$\class_{\nc}(G)$ has a natural metric on it: $\langle f, g \rangle
= \frac{1}{|G|} \sum_{\sigma \in G} f(\sigma)g(\sigma^{-1}).$  Let
$\Phi$ denote the standard additive isomorphism $$\Phi:
\class_{\nc}(G) \rTo Z\nc[G]$$ defined by linearly extending
$\Phi(e_{[\![\gamma]\!]}) = \sum_{\alpha \in [\![\gamma]\!]}
\alpha$.  The standard group-algebra product $\cdot$ in $Z\nc[G]$
pulls back via $\Phi$ to convolution $\star$ of functions.  And
the pushforward $\Phi_*\langle \ ,\ \rangle$ to $Z\nc[G]$ of the
metric  on $\class_{\nc}(G) $ agrees with the metric $(\ ,\ )$ on
$\nc[G]$ defined as $(\alpha,\beta) =
\frac{1}{|G|}\delta_{\alpha,\beta^{-1}}$, when $\alpha$ and
$\beta$ are in $G$.  With respect to this metric, $Z\nc[G]$ is a
Frobenius algebra.

\begin{crl}\label{crl:frob-iso}

The homomorphisms $$ (\ch,\ast,\eta) \rTo (\class_{\nc}(G),\star,
\langle \ ,\  \rangle ) \rTo^{\Phi} (Z\nc[G], \cdot, (\ ,\  )\ )$$
are isomorphisms of Frobenius algebras.
\end{crl}

\begin{proof}
It is well-known, and clear, that these maps are isomorphisms of
vector spaces.  The rest is a straightforward computation using
the definitions $$\eta_{[\![\sigma]\!][\![\beta]\!]} =
\langle\tau_0(e_{[\![\sigma]\!]}) \tau_0(e_{[\![\beta]\!]})
\tau_0(e_{[\![1]\!]})\rangle^G_0$$ and $$e_{[\![\alpha]\!]} \ast
e_{[\![\beta]\!]} := \sum_{[\![\sigma]\!],[\![\gamma]\!]}
\langle\tau_0(e_{[\![\alpha]\!]})\tau_0(e_{[\![\beta]\!]})
\tau_0(e_{[\![\sigma]\!]})\rangle^G_0
\eta^{[\![\sigma]\!][\![\gamma]\!]}e_{[\![\gamma]\!]}.$$
\end{proof}

These definitions are essentially the same as those given in
\cite{AGV}, but they differ \emph{a priori} from those of Chen and
Ruan \cite{CR1,Ru} for the ``orbifold Poincar\'{e} pairing'' and
``orbifold cup product.''  Nevertheless, it is easy to check that
both the geometry and the final calculations are identical to
those in \cite{CR1,Ru}.  For example, the three-fold multisectors
$X_{[\![\bga]\!]}$ of \cite[Defn 3.1.3]{Ru} are, in our case,
precisely the components $\M_{0,3}(\bg,[\![\bga]\!])$.

We are interested now in the corresponding cohomological field
theory (\cft) and the large phase space (when $a_i>0)$.

\begin{prop}\label{prop:cors}
The correlators $\langle\tau_{a_1}(e_{[\![\sigma_1]\!]}) \dots
\tau_{a_n}(e_{[\![\sigma_n]\!]})\rangle^G_g$ are related to the
usual correlators $\langle\tau_{a_1} \dots \tau_{a_n}\rangle_g$
(corresponding to the case of $G=\{ 1 \}$) by
$$\langle\tau_{a_1}(e_{[\![\sigma_1]\!]}) \dots
\tau_{a_n}(e_{[\![\sigma_n]\!]})\rangle^G_g = \langle\tau_{a_i}
\dots \tau_{a_n}\rangle_g \Omega^G_g (\bga)$$ where
        $$\Omega^G_{g}(\bga) =\frac{|\cx^G_{g} (\bga)|}{|G|},$$ and
$$\cx^G_{g} (\bga):=\{(\alpha_1, \dots, \alpha_g, \beta_1, \dots,
\beta_g, \sigma_1, \dots, \sigma_n) |\textstyle
\prod^g_{i=1}[\alpha_i,\beta_i] =\prod^n_{j=1}\sigma_j, \ \sigma_j
\in [\![\gamma_j]\!]\text{ for all } j \}.$$
\end{prop}

\begin{proof}
The sublocus
of $\mgnbar(\bg)$ where
$\prod^n_{i=1}ev^*_i(e_{[\![\gamma_i]\!]})$ is non-zero is
$\mgnbar(\bg,[\![\gamma_1]\!], \dots, [\![\gamma_n]\!]).$
Proposition \ref{point} shows that the degree of the forgetful map
$\pi:\mgnbar(\bg,[\![\gamma_1]\!], \dots, [\![\gamma_n]\!]) \rTo
\mgnbar$ is exactly $\Omega^G_g(\mathbf{\bga})$.  The proof now
follows from the projection formula, since $\psi_i$ on
$\mgnbar(\bg)$ is just the pullback $\pi^*\psi_i$ of the
corresponding class on $\mgnbar$.
\end{proof}

\begin{lm}\label{prop:omega}
The numbers $\Omega^G_g(\bga)$ depend only on the conjugacy
classes $[\![\gamma_i]\!]$, are independent of the ordering of the
$\gamma_i$ in $\bga$, and satisfy the following relations:
\begin{enumerate}
\item  \textbf{Cutting trees:}  For $g=g_1+g_2$ and $I \coprod J=
\{1,\dots,n\}$, let $\bga_I = (\gamma_{i_1}, \dots,
\gamma_{i_{|I|}})$ and $\bga_J=(\gamma_{j_1}, \dots,
\gamma_{j_{|J|}})$ $$\Omega^G_g (\bga) =
\Omega^G_{g_1}(\bga_I,\zeta) \eta^{[\![\zeta]\!][\![\xi]\!]}
\Omega^G_{g_2}(\xi,\bga_J)$$
\item \textbf{Cutting loops:} $$\Omega^G_g(\bga) =
\eta^{[\![\zeta]\!](\xi]\!]} \Omega^G_{g-1}(\zeta,\xi,\bga)$$
\item \textbf{Forgetting tails:} $$\Omega^G_g (\bga) = \Omega^G_g
(1,\bga)$$
\end{enumerate}

\end{lm}

\begin{proof}
Independence of conjugacy class representative and of order are
immediate from the definition, as is relation 3 (Forgetting
tails).  To prove relation 2, note that we may assume that
$[\![\zeta]\!] =[\![\xi^{-1}]\!].$  Let $\cx^G_{g} (\bga)$ be as
in Proposition~\ref{prop:cors}  and let $\cy$ be the set
\begin{align*}\cy:=&\{(\alpha'_1, \dots, \alpha'_{g-1}, \beta'_1,
\dots, \beta'_{g-1}, \sigma'_1,\dots, \sigma'_n, \sigma'_{n+1},
\sigma'_{n+2}) | \\ & \textstyle
\prod^{g-1}_{i=1}[\alpha'_i,\beta'_i]=\prod^{n+2}_{j=1}\sigma'_i,
\ \sigma'_j \in [\![\gamma_j]\!] \text{ for } j \leq n,\
\sigma'_{n+2} \in [\![{\sigma'}^{-1}_{n+1}]\!]\}.\end{align*}
Define a map $f:\cx_g^G(\bga) \rTo \cy,$ taking
$(\alpha_1, \dots, \alpha_g, \beta_1, \dots, \beta_g, \sigma_1,
\dots, \sigma_n)$ in $\cx^G_{g} (\bga)$  to $(\alpha_1, \dots,
\alpha_{g-1}, \beta_1, \dots, \beta_{g-1}, \sigma_1, \dots,
\sigma_n, \alpha_g\beta_g\alpha^{-1}_g,\beta^{-1}_g)$ in $\cy$.

For a given conjugacy class $[\![\psi]\!]$ in $G$, the map $f$,
restricted to the subset of $\cx^G_{g} (\bga)$ where $\beta_g \in
[\![\psi]\!]$, takes  $|C(\psi)|$ elements of $\cx^G_{g} (\bga)$
to one element of $\cy$. Moreover, given $\zeta,\xi \in G$ such
that $[\![\zeta^{-1}]\!] = [\![\xi]\!]$, if we let $\beta_g =
\xi^{-1}$ and $\alpha_g$ be an element such that $\alpha_g \beta_g
\alpha^{-1}_g= \zeta$, then it is clear that the map $f$ is
surjective. This shows that $\Omega^G_g (\bga) =
\sum_{[\![\zeta]\!]} |C(\zeta)|
\Omega^G_{g-1}(\bga,\zeta,\zeta^{-1})$, which proves (2). The proof
of (1) is a similar, straightforward argument.
\end{proof}

\begin{thm}
Let $\Lambda^G_{g,n} :\ch^{\otimes n} \rTo H^*(\mgnbar)$ be
defined as $\Lambda^G_{g,n}(e_{[\![\gamma_1]\!]}\otimes \dots
\otimes e_{[\![\gamma_n]\!]})=
\pi_*(ev^*_1(e_{[\![\gamma_1]\!]})\dots
ev^*_n(e_{[\![\gamma_n]\!]}))=\Omega^G_g(\bga).$ The collection
$(\ch,\eta,\Lambda^G,e_{[\![1]\!]})$ is a \cft\ with flat
identity.
\end{thm}

\begin{proof}
It is straightforward to see that
$\Lambda^G_{g,n}(e_{[\![\gamma_1]\!]}\otimes \dots \otimes
e_{[\![\gamma_n]\!]}) =\Omega^G_n (\bga)$.  The \cft\ axioms
follow immediately from Lemma ~\ref{prop:omega}.
\end{proof}

The \cft\ axioms also follow from the cutting and forgetting tails
axioms that hold for the virtual fundamental class \cite{CR2},
which is trivial for $\mgnbar(\bg)$.  Alternatively, the axioms can
be seen directly from the geometry by carefully accounting for
ramification of $\pi$ on the boundary and accounting for the
degree of the morphism $$\mgnbar(\bg) \times_{\mgnbar}
(\M_{g,n_1+1} \textstyle \coprod \M_{g_2,n_2+1}) \rTo
\M_{g_1,n_1+1} (\bg) \coprod \M_{g_2,n_2+1} (\bg)$$ as in
\cite[\S4.2]{JKV}.

\subsection{K-Theory} The orbifold K-theory of $\bg$ is simply
the representation ring $R_G$, since a vector bundle on $\bg$ is a
vector space with linear $G$-action.  The Chern character (see
\cite{AR}) $ch:K^*_{orb} \otimes \nc \rTo H^*_{orb} (\bg,\nc)
=\ch$ is easily seen to be the composite $$R_G \otimes \nc
\rTo^\chi \class_{\nc}(G) \cong \ch$$ of the trace map $\chi$ and
the obvious {additive} isomorphism from class functions to $\ch$.

The trace map $\chi$ is a ring isomorphism, but as stated in
Corollary~\ref{crl:frob-iso}, the orbifold product $\ast$ on $\ch$
corresponds to convolution $\star$ of class functions rather than
multiplication.  Thus the Chern character is only an additive
isomorphism.

\subsection{Functoriality}

\subsubsection{Morphisms}\mbox{}\\
It is interesting to note that K-theory has functoriality
properties that $H^*_{orb}$ does not enjoy.  In particular, a
homomorphism of groups $G \rTo H$ gives a morphism $\bg \rTo \cb
H$ which induces a ring homomorphism $$K^*_{orb}(\cb H) \rTo
K^*_{orb}(\bg),$$ corresponding to the obvious homomorphism of
representation rings. But the center $Z \nc [H]$ generally has no
ring homomorphism to $Z \nc[G]$.

\subsubsection{Tensor products}\mbox{}\\
Tensor products arise in this theory in
at least
two ways.

\begin{itemize}
\item \textbf{$\cb (G\cross H)$}

For any two finite groups $G$ and $H$, the classifying stack $\cb (G\cross H)$
splits up as a product
$$\cb (G\cross H) =  \bg \cross \cb H,$$ and although the moduli
stack of the product is not quite the product of the moduli stacks
of the factors, it is easy to see that the corresponding \cft\ is
the tensor product of the two components. That is, a
straightforward check shows that if $\bga =
(\gamma_1,\dots,\gamma_n) \in G^n$ and $\bsig =
(\sigma_1,\dots,\sigma_n) \in H^n$ then  $$\Omega^{G\cross H}_g(
(\gamma_1,\sigma_1), \dots, (\gamma_n,\sigma_n)) =
\Omega^G_g(\bga) \Omega^H_g(\bsig).$$

\item \textbf{$[X/G]$ with trivial $G$ action}

If $X$ is a smooth projective variety with trivial $G$ action,
then the quotient stack $[X/G]$ is isomorphic to the product $X
\cross \bg$.

\begin{prop}
For a smooth projective variety $X$ with finite group $G$ acting
trivially,  the \cft\ arising from stable maps into the orbifold
$[X/G]$ is simply the tensor product of the \cft\ arising from
stable maps into $X$ and the \cft\ arising from stable maps into
$\bg$.
\end{prop}
\begin{proof}
This follows from the fact that the degree of the forgetful maps
$\M_{g,n}(X\times \bg,[\![\gamma_1]\!], \dots,
[\![\gamma_n]\!])\to \M_{g,n}(X)$ and
$\M_{g,n}(\bg,[\![\gamma_1]\!], \dots, [\![\gamma_n]\!])\to
\M_{g,n}$ are both equal to $\Omega_g(\bga)$.
\end{proof}

\end{itemize}

\section{Semisimplicity and Virasoro Algebras}

In this section, the summation convention is NOT used on any
subscripts or superscripts $\alpha$ or $\alpha_i$, although it is
applied to all other variables.

\subsection{Semisimple Frobenius algebras} Let $V$ be any $r$-dimensional
Frobenius algebra with multiplication $*$, metric $\eta$, and
identity element $\mathbf{1}$. It is said to be a \emph{semisimple
Frobenius algebra} if there exists a \emph{canonical basis}
$\{\,f_\alpha\,\}_{\alpha=1}^r$ such that for all $\alpha_1,
\alpha_2= 1,\ldots,r$,
\begin{equation}
\label{ss1} f_{\alpha_1} * f_{\alpha_2} =
\delta_{\alpha_1,\alpha_2} f_{\alpha_1}
\end{equation}
and
\begin{equation}
\label{ss2} \eta(f_{\alpha_1},f_{\alpha_2}) =
\delta_{{\alpha_1},{\alpha_2}} \nu_{\alpha_1}
\end{equation}
for some nonzero numbers $\nu_\alpha$. The identity element
satisfies
\[
\mathbf{1} = \sum_{\alpha=1}^r f_\alpha.
\]
As discussed before, the Frobenius algebra
$(\ch,\eta,*,e_{[\![1]\!])})$ can be identified with $Z \nc[G]$
and the latter is a semisimple Frobenius algebra with canonical
basis given as follows.

\begin{prop}
Let $\{\,V_\alpha\,\}_{\alpha=1}^r$ be the set of irreducible
representations of $G$ and let $\chi_\alpha$ denote the character
of $V_\alpha$.  For all $\alpha=1,\ldots,r$, the elements
\[
f_\alpha := \frac{\dim V_\alpha}{|G|}\sum_{g\in G}
\chi_{\alpha}(g^{-1}) g
\]
form a basis of $Z\nc[G]$ and satisfy equations (\ref{ss1}) and
(\ref{ss2}), where for all $\alpha=1, \ldots, r$,
\[
\nu_\alpha = \left(\frac{\dim V_\alpha}{|G|}\right)^2.
\]
\end{prop}
\begin{proof}

It is clear from the definition that the $f_{\alpha}$ lie in
$Z\nc[G]$.  The fact that they satisfy (\ref{ss1}) follows from
\cite[\S2.4]{FH}. That equation (\ref{ss2}) holds for the given
values of $\nu_{\alpha}$ is a straightforward computation.
\end{proof}

The results in the remainder of this section hold for any
semi-simple Frobenius algebra, since any Frobenius algebra is a
\cft.

\subsection{The Potential Function}

We will now calculate the correlators for our theory in the
canonical coordinates.

\begin{prop}\label{correlators}
Let $(g,n)$ be any stable pair, where $n\geq 1$. If $\alpha_i =
\alpha$ for all $i=1,\ldots,n$ then for all $a_1,\ldots, a_n\geq
0$, we have
\begin{equation}
\corG{\tau_{a_1}(f_{\alpha})\cdots\tau_{a_n}(f_{\alpha})}_g =
\nu_\alpha^{1-g} \cor{\tau_{a_1}\cdots\tau_{a_n}}_g;
\end{equation}
otherwise, we have
$\corG{\tau_{a_1}(f_{\alpha_1})\cdots\tau_{a_n}(f_{\alpha_n})}_g =
0$. Furthermore, when $n=0$, we have $\corG{}_g = 0$.
\end{prop}
\begin{proof}
Of course, $\corG{}_g = 0$ holds for dimensional reasons.

The proof for the rest of the proposition follows by degenerating
to curves whose irreducible  components are all three-pointed,
genus-zero curves, where the proposition is easily verified, and
then calculating the general correlators from the cutting axioms
for \cfts.

More explicitly, the correlator is simply
$$\corG{\tau_{a_1}(f_{\alpha_1})\cdots\tau_{a_n}(f_{\alpha_n})}_g
= \Lambda^G_{g,n}(f_{\alpha_1}\tensor\dots\tensor f_{\alpha_n})
\cor{\tau_{a_1}\cdots\tau_{a_n}}_g,$$ where $(\ch,\eta,\Lambda^G)$
is our \cft.  The definition of $\ast$ gives
\begin{align*}
\Lambda^G_{0,3}(f_{\alpha_1}\tensor f_{\alpha_2}\tensor
f_{\alpha_3}) &= \eta(f_{\alpha_1}\ast f_{\alpha_2},
f_{\alpha_3})\\ &=
\delta_{\alpha_1,\alpha_2}\delta_{\alpha_1,\alpha_3}\nu_{\alpha_1}.
\end{align*}

Now proceed by induction on the genus and number of marked points.
Each application of the \emph{cutting trees} axiom
(for a $3$-pointed, genus-zero vertex)
leaves the genus unchanged, and reduces the number of marked
points by one, but contributes nothing to the final result, since
the node (cut edge) contributes the inverse metric---a factor of
$\nu_{\alpha}^{-1}$---and the $3$-pointed, genus-zero, irreducible
component (vertex of the dual graph) contributes  a factor of
$\nu_{\alpha}$.

Each application of the \emph{cutting loops} axiom increases the
number of marked points by $2$, reduces the genus by $1$, and
contributes the inverse of the metric---namely,
$\nu_{\alpha}^{-1}$---to the final result.

\end{proof}

The \emph{large phase space potential} is defined by $\PhiG(\bt) =
\sum_g \PhiG_g(\bt) \lambda^{2g-2}$ in
$\lambda^{-2}\nc[[\bt,\lambda]]$, where $\PhiG_g(\bt) :=
\corG{\exp(\btau\cdot\bt)}_g$, where $\bt\cdot\btau = \sum_{a,m}
t_a^m \tau_a(h_m)$, and where $\{ h_m \}$ is any basis for $\ch$.
Let $\ZG := \exp(\PhiG)$.

When $G=\{\mathbf{1}\}$ is the trivial group, we denote by $\Phi
:= \Phi^{\{\mathbf{1}\}}$, the potential of the Gromov-Witten
invariants of a point. Similarly, we let $Z :=
Z^{\{\mathbf{1}\}}$.

\begin{prop}\label{factorization}
Let $\bu$ be formal variables $\{\,u_a^\alpha\,\}$ for all integers $a\geq
0$ and $\alpha=1,\ldots,r$ associated to the canonical basis
$\{\,f_\alpha\,\}$. For each $\alpha = 1,\ldots,r$, let
$\tbu^\alpha$ be formal variables $\{\,\tu_a^\alpha\,\}$, where
$a\geq 0$ and $\tu_a^\alpha := (\nu_\alpha)^{\frac{1-a}{3}}
u_a^\alpha$. Then we have
\begin{equation}
\PhiG(\bu) = \sum_{\alpha=1}^r \Phi(\tbu^\alpha).
\end{equation}
\end{prop}
\begin{proof}
This follows from Proposition \ref{correlators} and dimensional
considerations.
\end{proof}

\subsection{Virasoro algebras}

For each $\alpha=1,\ldots,r$ and $n\geq -1$, let
\begin{eqnarray*}
L_n^{(\alpha)} &:=& -\frac{(2n+3)!!}{2^{n+1}}
\frac{\partial}{\partial \tu^\alpha_{n+1}} + \sum_{i=0}^\infty
\frac{(2i+2n+1)!!}{(2i-1)!! 2^{n+1}}
\tu_i^\alpha\frac{\partial}{\partial \tu_{i+n}^\alpha} \\ &+&
\frac{\lambda^2}{2}\sum_{i=0}^{n-1} \frac{(2i+1)!!
(2n-(2i+1))!!}{2^{n+1}} \frac{\partial^2}{\partial \tu_i^\alpha
\partial \tu_{n-1-i}^\alpha} \\ &+& \delta_{n,-1}
\frac{\lambda^{-2}}{2} \tu_0^\alpha \tu_0^\alpha + \delta_{n,0}
\frac{1}{16},
\end{eqnarray*}
where$(2n-1)!! := 1 \cdot 3\cdot 5\ldots\cdot (2n-1)$. These
operators satisfy
\[
[L_m^{(\alpha)},L_n^{(\beta)}] = (m-n) L^{(\alpha)}_{m+n}
\delta_{\alpha,\beta}
\]
for all $m,n\geq -1$ and $\alpha,\beta \in \{1,\ldots,r\}$,
forming $r$ commuting copies of ``half'' of the Virasoro algebra.

Moreover, if $\{\,b_m\,\}$ is any basis for $\ch$, where
$m=0,\ldots,r-1$ such that $b_0 = \mathbf{1}$, and if $\bt =
\{\,t_a^m\,\}$ are the associated formal parameters, then there
are also operators for all $n\geq -1$ given by
\begin{eqnarray*}
L_n &:=& -\frac{(2n+3)!!}{2^{n+1}} \frac{\partial}{\partial
t^0_{n+1}} + \sum_{i=0}^\infty \frac{(2i+2n+1)!!}{(2i-1)!!
2^{n+1}} \left(\sum_{m} t_i^m\frac{\partial}{\partial
t^m_{i+n}}\right) \\ &+& \frac{\lambda^2}{2}\sum_{i=0}^{n-1}
\frac{(2i+1)!! (2n-(2i+1))!!}{2^{n+1}}
\left(\sum_{m_1,m_2}\eta^{m_1 m_2}\frac{\partial^2}{\partial
t_i^{m_1} \partial t_{n-1-i}^{m_2}}\right)
\\ &+&
\delta_{n,-1} \frac{\lambda^{-2}}{2} \left(\sum_{m_1,m_2}
\eta_{m_1 m_2} t_0^{m_1} t_0^{m_2}\right) + \delta_{n,0}
\frac{r}{16},
\end{eqnarray*}
satisfying $[L_k,L_n] = (k-n) L_{k+n}$ for any $k,n\ge -1$.

\begin{prop}
For all $\alpha \in \{ 1,\ldots,r\}$ and $n\geq -1$,
\begin{equation}\label{fullvirasoro}
L_n^{(\alpha)} \ZG = 0.
\end{equation}
These equations completely determine $\ZG$. Furthermore, for all
$n\geq -1$,
\begin{equation}\label{diagonalvirasoro}
L_n \ZG = 0.
\end{equation}
\end{prop}
\begin{proof}
Equation (\ref{fullvirasoro}) follows from Proposition
\ref{factorization} and the Kontsevich-Witten theorem \cite{Ko,Wi}
for the case of $G=\{\mathbf{1}\}$. Equation
(\ref{diagonalvirasoro}) follows from (\ref{fullvirasoro}) and the
identity
\[
L_m = \sum_{\alpha=1}^r (\nu_\alpha)^{-\frac{m}{3}}
L_m^{(\alpha)}.
\]
\end{proof}

\begin{rem} Equation \ref{diagonalvirasoro}
is a verification of the Virasoro conjecture for $\bg$ \cite{EHX}
and can also be regarded as an example of \cite{Gi}.
\end{rem}

\subsection{KdV hierarchies}

Let $\ccor{A}^G_g := \ccor{A\exp(\mathbf{t}\cdot\btau)}^G_g$ and
$\ccor{A}^G := \sum_g \ccor{A\exp(\mathbf{t}\cdot\btau)}_g^G
\lambda^{2g-2}.$ The superscript $G$ will be suppressed when
$G=\{\mathbf{1}\}$, the trivial group.

\begin{prop}
Let $\{\,e_1,\ldots,e_m\,\}$ be any basis for $\ch$. For all $v$
in $\ch$ and $a\geq 0$, the following equation holds:
\begin{equation} \label{KdV}
\begin{split}
 &(2a+1) \lambda^{-2}  \ccor{\tau_a(v)\tau_0(e_{m_1})\tau_0(e_{m_2})}^G
\eta^{m_1 m_2} =  \\ & \ccor{\tau_{a-1}(v) \tau_0(e_{m_1})}^G
\eta^{m_1 m_2} \ccor{\tau_0(e_{m_2})
\tau_0(e_{m_3})\tau_0(e_{m_4})}^G \eta^{m_3 m_4}+ \\ &2
\ccor{\tau_{a-1}(v)\tau_0(e_{m_1})\tau_0(e_{m_3})}^G\eta^{m_1 m_2}
\eta^{m_3 m_4} \ccor{\tau_0(e_{m_2}),\tau_0(e_{m_4})}^G \eta^{m^3
m^4} +\\ & \frac{1}{4} \ccor{\tau_{a-1}(v) \tau_0(e_{m_1})
\tau_0(e_{m_2}) \tau_0(e_{m_3}) \tau_0(e_{m_4})}^G \eta^{m_1 m_2}
\eta^{m_3 m_4}.\\
\end{split}
\end{equation}
Equation (\ref{KdV}) and the fact that $L_{-1}^{(\alpha)} Z^G = 0$
for all $\alpha=1,\ldots,r$ completely determine $\Phi^G$.
\end{prop}

\begin{proof}
When $G=\{\mathbf{1}\}$, the trivial group, equation (\ref{KdV})
reduces to
\begin{equation}
\begin{split}
 &(2a+1) \lambda^{-2}  \ccor{\tau_a\tau_0\tau_0} = \\
 &  \ccor{\tau_{a-1} \tau_0}
\ccor{\tau_0 \tau_0\tau_0} + 2 \ccor{\tau_{a-1}\tau_0\tau_0}
\ccor{\tau_0,\tau_0}  + \frac{1}{4} \ccor{\tau_{a-1} \tau_0 \tau_0
\tau_0 \tau_0}. \\
\end{split}
\end{equation}
This equation is the Kontsevich-Witten theorem \cite{Ko,Wi}.
Witten also showed that this equation together with $L_{-1} Z = 0$
completely determines $\Phi$.

To prove the formula for general $G$, choose a canonical basis
$\{\,f_\alpha\,\}_{\alpha=1}^r$ and let $v = f_\alpha$. Consider
the terms of equation (\ref{KdV}) proportional to
$\lambda^{2g-4}$. By Proposition (\ref{correlators}), one obtains
the terms of equation (\ref{KdV}) proportional to $\lambda^{2g-4}$
up to an overall scalar factor.
\end{proof}

\begin{rem}
Equation \ref{KdV} is a simultaneous solution of $r$ commuting KdV
hierarchies with time parameters $\widetilde{u}_a^\alpha$ for
$a\geq 0$ and $\alpha=1,\ldots,r$.
\end{rem}

\bibliographystyle{amsplain}

\providecommand{\bysame}{\leavevmode\hbox
to3em{\hrulefill}\thinspace}

\end{document}